# Robust Fixed-order Dynamic Output Feedback Controller Design for Fractional-order Systems

Pouya Badri[1], Mahdi Sojoodi[1*]

[1]Advanced Control Systems Laboratory, School of Electrical and Computer Engineering, Tarbiat Modares University, Tehran, Iran.
[*]sojoodi@modares.ac.ir

**Abstract:** This paper deals with designing a robust fixed-order dynamic output feedback controller for uncertain fractional-order linear time invariant (FO-LTI) systems by means of linear matrix inequalities (LMIs). Our purpose is to design a low-order controller that stabilizes the fractional-order linear system in the presence of model uncertainties. No limiting constraint on the state space matrices of the uncertain system is assumed in the design procedure. Furthermore, adopting the most complete model of linear controller, with direct feedthrough parameter, does not disturb the LMI-based approach of developing robust stabilizing control. Eventually, the authors illustrate the advantages of the proposed method by some examples and their numerical simulation.

## 1. Introduction

In recent years great effort has been devoted to the study of applications of fractional calculus in the modelling and control of various real-world systems [1-3]. Fractional-order models can more precisely describe systems that have responses with long memory transients. Moreover, it has long been recognized that many natural, biological, and practical engineering systems have inherent attributes that can be better explained by fractional-order models [4-6]. Therefore, controller-designing methods for systems modelled by fractional-order dynamics are among the most widely discussed issues in the literature. Using fractional-order controllers is an appropriate procedure for controlling such systems due to their high performance and robustness [7,8]. Since, in the most of control systems, stability is the first objective to be accomplished, and due to uncertain models caused by neglected dynamics, parametrical variations in time, uncertain physical parameters, and so on, robust stability and stabilization have become a fundamental issue for all control systems including fractional-order systems [9,10].

Robust stability and stabilization of fractional-order control systems are investigated in [9-13]. Necessary and sufficient conditions for the stability and stabilization of fractional-order interval systems are presented in [9]. In addition, necessary and sufficient conditions for the robust stability of general interval fractional-order system are investigated in [14], where interval uncertainties exist both in the coefficients and orders of the fractional-order system. The existence conditions and design procedures of the static state feedback controller, static output feedback controller and observer-based controller for asymptotically stabilizing fractional-order linear systems with positive real uncertainty are presented in [10] where designing a static output feedback controller is possible with the constraint on the output matrix of the uncertain system to be of full-row rank. In [15] satisfaction of some frequency domain specifications in tuning fractional-order proportional-derivative controllers guarantees robustness of the closed-loop system against DC gain variations.

In [16] the problem of the observer-based robust control for a class of Lipschitz nonlinear fractional-order systems is investigated. The proposed method therein, is used for stabilization of fractional-order interval systems with fractional order $0 < \alpha < 1$ by assuming that the output matrix of the uncertain system is of full row rank. Moreover, the concentration of [17] is on the stabilization of fractional-order systems subject to bounded uncertainties where the system uncertainties are randomly distributed in the state matrix $A$ and the output matrix $C$ and fractional differentiation order is between zero and one. In this paper, in order to solve the problem of observer-based stabilization of fractional-order uncertain system, it is assumed that all individual possible pairs of $(A + \Delta A, C + \Delta C)$ is observable in the sense of Kalman. Furthermore, in [18] by using singular value decomposition and LMI techniques, robust control of fractional-order interval linear systems were considered. The proposed method therein, is used for stabilization of fractional-order interval systems with fractional order $1 \leq \alpha < 2$ by assuming that the output matrix of the uncertain system is of full row rank.

In the most of mentioned works, state feedback controller is employed which entails having all individual states. On the other hand, in some cases measuring all states seems impossible due to economic issues or physical constraints [19], where using output feedback control seems to be helpful. It is important to note that dynamic feedback controllers are preferable to static ones due to their more effective control performances alongside with more degrees of freedom in achieving control objectives [20].

A great number of the controller design methods lead to high order controllers suffering from expensive implementation, undesirable reliability, high fragility, maintenance difficulties, and potential numerical errors. Since plant or controller order reduction methods do not always guarantee the closed-loop performance, a possible solution is to design a controller with a low and fixed order, chosen by the designer to be as small as possible according to the requirements of the problem [19].

It is notable that using the interval uncertainty description, polytopic uncertainty description, and norm-bounded uncertainty description one can only capture gain uncertainty and in the case that phase uncertainty information is available, these uncertainty descriptions may lead to conservative results [10,21,22]. A suitable way for considering phase information is to use the positivity theorem or, in other words, to model the uncertainty through a positive real uncertainty matrix [10, 21–24]. It should be noted that positive real uncertainty exists in many real systems, and the robust stability and stabilization problem of integer-order systems with positive real uncertainty has been investigated in [10, 21–24].

As far as we know, there is few reported work on the analytically designing of a stabilizing dynamic output feedback controller for fractional-order systems with positive real uncertainty in the literature. In this paper, the robust stabilization of fractional-order linear systems with positive real uncertainty through a dynamic output feedback controller with a predetermined order is investigated. In spite of the complexity of assuming the most complete model of linear controller with direct feedthrough parameter, LMI approach of developing robust control is preserved, which is convenient to be used in practice due to several effective convex optimization Parsers and Solvers that can be used to assess the feasibility of the LMI constraints and obtain design parameters. Moreover, designing an output feedback controller does not impose any limiting constraint on the state space matrices of the uncertain systems which was considered in [10].

The remainder of this paper is organized as follows. In section 2, some preliminaries and the problem formulation are presented. Robust stabilizing conditions of uncertain fractional-order systems with positive real uncertainty via fixed-order dynamic output feedback controller alongside with the design procedures of the corresponding controllers are derived in Section 3. Some numerical examples are provided in Section 4 to illustrate the effectiveness of the proposed method. Eventually, section 5 concludes the paper.

**Notations**: In this paper $A \otimes B$ denotes the kronecker product of matrices $A$ and $B$ and by $M^T$, $\bar{M}$ and $M$, we denote the transpose, the conjugate, and the transpose conjugate of $M$, resepectively. $\bar{z}$ represents the conjugate of the scalar number $z$ and $Sym(M)$ stands for $M + M^*$. The notation ● is the symmetric component in matrix and ↑ is the symbol of pseudo inverse. Matrix $I$ denotes the identity matrix with appropriate dimensions and $i$ denotes the imaginary unit.

## 2. Preliminaries and problem formulation

In this section, some basic concepts and lemmas of fractional-order calculus together with positive real uncertainty are presented.

Consider the following uncertain fractional-order linear time-invariant (FO-LTI) system:

$$D^\alpha x(t) = (A + \Delta A)x(t) + (B + \Delta B)u(t), 0 < \alpha < 2$$

$$y(t) = Cx(t) \qquad (1)$$

in which $x \in R^n$ denotes the pseudo-state vector, $u \in R^l$ is the control input, and $y \in R^m$ is the output vector. Furthermore $A \in R^{n \times n}, B \in R^{n \times l}$, and $C \in R^{m \times n}$ are known constant matrices with appropriate dimensions. $\Delta A$ and $\Delta B$ are time-invariant matrices with parametric uncertainty of the following form [22,23]

$$[\Delta A \quad \Delta B] = M\Delta(\zeta)[N_1 \quad N_2], \qquad (2)$$

$$\Delta(\zeta) = F(\zeta)[I + JF(\zeta)]^{-1}, \qquad (3)$$

$$sym(J) > 0, \qquad (4)$$

where $M \in R^{n \times m_0}$, $N_1 \in R^{m_0 \times n}$, $N_2 \in R^{m_0 \times l}$ and $J \in R^{m_0 \times m_0}$ are real known matrices. The uncertain matrix $F(\zeta) \in R^{m_0 \times m_0}$ satisfies

$$Sym\{F(\zeta)\} \geq 0, \qquad (5)$$

where $\zeta \in \Omega$ with $\Omega$ being a compact set.

**Remark 1**: Condition (4) guarantees that $I + JF(\zeta)$ is invertible for all $F(\zeta)$ satisfying (5). Therefore $\Delta(\zeta)$ in (3) is well defined [10].

In this paper, we adopt the following Caputo definition for fractional derivatives of order α of function $f(t)$ since initial values of classical integer-order derivatives with explicit physical concepts are utilizable using the Laplace transform of the Caputo derivative [25]:

$$^C_aD^\alpha_t f(t) = \frac{1}{\Gamma(m-\alpha)} \int_a^t (t-\tau)^{m-a-1} \left(\frac{d}{d\tau}\right)^m f(\tau)d\tau$$

where $\Gamma(\cdot)$ is Gamma function defined by $\Gamma(\epsilon) = \int_0^\infty e^{-t}t^{\epsilon-1}dt$ and $m$ is the smallest integer that is equal to or greater than $\alpha$.

In order to study the stability of fractional-order systems and obtain main results the following lemmas are required.

**Lemma 1** [26]: Let $A \in R^{n \times n}$, $0 < \alpha < 1$ and $\theta = \alpha\pi/2$. The fractional-order system $D^\alpha x(t) = Ax(t)$ is asymptotically stable if and only if there exists a positive definite Hermitian matrix $X = X^* > 0, X \in C^{n \times n}$ such that

$$(rX + \bar{r}\bar{X})^T A^T + A(rX + \bar{r}\bar{X}) < 0, \qquad (6)$$

where $r = e^{\theta i}$.

**Lemma 2** [27]: Let $A \in \mathcal{R}^{n \times n}$, $1 \leq \alpha < 2$ and $\theta = \pi - \alpha\pi/2$. The fractional-order system $D^\alpha x(t) = Ax(t)$ is asymptotically stable if and only if there exists a positive definite matrix $X \in \mathcal{R}^{n \times n}$ such that

$$\begin{bmatrix} (A^TX + XA)\sin\theta & (XA - A^TX)\cos\theta \\ \bullet & (A^TX + XA)\sin\theta \end{bmatrix}, \qquad (7)$$

defining

$$\theta = \begin{bmatrix} \sin\theta & -\cos\theta \\ \cos\theta & \sin\theta \end{bmatrix}, \qquad (8)$$

and with this in mind that $A$ is similar to $A^T$, inequality (7) can be expressed as follows

$$Sym\{\Theta \otimes (AX)\} < 0. \quad (9)$$

**Lemma 3** [22]: Let $\Omega = \{\Delta \in \mathcal{R}^{m_0 \times m_0} | \Delta \text{ is subjected to } (3) - (5)\}$. Then

$$\Omega = \{\Delta(\zeta) \in R^{m_0 \times m_0} | det(I - \Delta(\zeta)J)) \neq 0 \text{ and } \Delta(\zeta)Sym\{J\}\Delta(\zeta)^T \leq Sym\{\Delta(\zeta)\}\}. \quad (10)$$

## 3. Main results

The main objective of this paper is to design a robust fixed-order dynamic output feedback controller that asymptotically stabilizes the uncertain FO-LTI system (1) in terms of linear matrix inequalities (LMIs). Therefore, the following dynamic output feedback controller is presented

$$D^\alpha x_C(t) = A_C x_C(t) + B_C y(t), \quad 0 < \alpha < 2$$

$$u(t) = C_C x_C(t) + D_C y(t) \quad (11)$$

with $x_C \in \mathcal{R}^{n_c}$, in which $n_c$ is the arbitrary order of the controller and $A_C$, $B_C$, $C_C$, and $D_C$ are corresponding matrices to be designed.

The resulted closed-loop augmented FO-LTI system using (1) and (11) is as follows

$$D^\alpha x_{Cl}(t) = A_{Cl} x_{Cl}(t), \quad 0 < \alpha < 2 \quad (12)$$

where

$$x_{Cl}(t) = \begin{bmatrix} x(t) \\ x_C(t) \end{bmatrix},$$

$$A_{Cl} = \begin{bmatrix} A + \Delta A + (B + \Delta B)D_C C & (B + \Delta B)C_C \\ B_C C & A_C \end{bmatrix} \quad (13)$$

**Theorem 1**: Considering closed-loop system in (12) with $0 < \alpha < 1$, and a positive definite Hermitian matrix $P = P^*$ in the form of

$$P = diag(P_S, P_C) \quad (14)$$

with $P_S \in \mathcal{C}^{n \times n}$ and $P_C \in \mathcal{C}^{n_c \times n_c}$ and a real scalar constant $\gamma > 0$ alongside with matrices $T_i, i = 1, ..., 4$ exist such that the following LMI constrain become feasible

$$\begin{bmatrix} \Pi_{11} & \widetilde{M} & \Pi_{13} \\ \bullet & -\gamma I & \gamma I \\ \bullet & \bullet & -Sym(J) - \gamma I \end{bmatrix} < 0, \quad (15)$$

in which

$$\Pi_{11} = \begin{bmatrix} p_{11} & p_{12} \\ p_{21} & p_{22} \end{bmatrix}, \Pi_{13} = \begin{bmatrix} q_1 \\ q_2 \end{bmatrix}$$

$$p_{11} = A(rP_S + \bar{r}\overline{P_S}) + (rP_S + \bar{r}\overline{P_S})^T A^T + BT_4 + T_4^T B^T,$$

$$p_{12} = BT_3 + T_2^T, p_{21} = T_2 + T_3^T B^T, p_{22} = T_1 + T_1^T,$$

$$q_1 = (rP_S + \bar{r}\overline{P_S})^T N_1^T + T_4^T N_2^T, q_2 = T_3^T N_2^T \quad (16)$$

where $\theta = (1 - \alpha)\pi/2$ and $r = e^{i\theta}$ then, the dynamic output feedback controller parameters of

$$A_C = T_1(rP_C + \bar{r}\overline{P_C})^{-1}, B_C = T_2(rP_S + \bar{r}\overline{P_S})^{-1}C^\uparrow,$$

$$C_C = T_3(rP_C + \bar{r}\overline{P_C})^{-1}, D_C = T_4(rP_S + \bar{r}\overline{P_S})^{-1}C^\uparrow, \quad (17)$$

make the closed-loop system in (12) asymptotically stable.

**Proof**: Let $A_{Cl} = A_{0Cl} + A_{\Delta Cl}$, with

$$A_{0Cl} = \begin{bmatrix} A + BD_C C & BC_C \\ B_C C & A_C \end{bmatrix},$$

$$A_{\Delta Cl} = \begin{bmatrix} \Delta A + \Delta BD_C C & \Delta BC_C \\ 0 & 0 \end{bmatrix} = \widetilde{M}\Delta\widetilde{N}$$

$$\widetilde{M} = \begin{bmatrix} M \\ 0 \end{bmatrix}, \widetilde{N} = [N_1 + N_2 D_C C \quad N_2 C_C]. \quad (18)$$

It follows from Lemma 1 that the uncertain fractional-order closed-loop system (12) with $0 < \alpha < 1$ is asymptotically stable if there exists a positive definite Hermitian matrix $X = X^*$, $X \in \mathcal{C}^{(n+n_c) \times (n+n_c)}$ in a way that

$$(rX + \bar{r}\bar{X})^T A_{Cl}^T + A_{Cl}(rX + \bar{r}\bar{X}) < 0 \Leftrightarrow$$

$$Sym\{A_{Cl}(rX + \bar{r}\bar{X})\} < 0 \Leftrightarrow$$

$$Sym\{(A_{0Cl} + A_{\Delta Cl})(rX + \bar{r}\bar{X})\} < 0 \Leftrightarrow$$

$$Sym\{A_{0Cl}(rX + \bar{r}\bar{X})\} + Sym\{\widetilde{M}\Delta(\zeta)\widetilde{N}(rX + \bar{r}\bar{X})\} < 0 \quad (19)$$

Introducing

$$E = Sym(J),$$

$$F = E^{-\frac{1}{2}}\left(\rho^{-1}\widetilde{M}^T + \rho\widetilde{N}(rX + \bar{r}\bar{X})\right) - E^{\frac{1}{2}}\Delta^T(\zeta)\rho^{-1}\widetilde{M}^T, \quad (20)$$

it follows from Lemma 3 that $Sym\{\Delta(\zeta)\} - \Delta(\zeta)E\Delta^T(\zeta) > 0$ and then the following inequality holds for any $\rho > 0$

$$-F^T F \leq 0 \Leftrightarrow -\left(E^{-\frac{1}{2}}\left(\rho^{-1}\widetilde{M}^T + \rho\widetilde{N}(rX + \bar{r}\bar{X})\right) - E^{\frac{1}{2}}\Delta^T(\zeta)\rho^{-1}\widetilde{M}^T\right)^T \left(E^{-\frac{1}{2}}\left(\rho^{-1}\widetilde{M}^T + \rho\widetilde{N}(rX + \bar{r}\bar{X})\right) - E^{\frac{1}{2}}\Delta^T(\zeta)\rho^{-1}\widetilde{M}^T\right) < 0 \Leftrightarrow$$

$$-Sym\{\widetilde{M}E^{-1}\widetilde{N}(rX + \bar{r}\bar{X})\} - \rho^{-2}\widetilde{M}E^{-1}\widetilde{M}^T - \rho^2(rX + \bar{r}\bar{X})^T \widetilde{N}^T E^{-1} N (rX + \bar{r}\bar{X}) + Sym\{\widetilde{M}\Delta(\zeta)\widetilde{N}(rX + \bar{r}\bar{X})\} + \rho^{-2}\widetilde{M}\left(Sym\{\Delta(\zeta)\} - \Delta(\zeta)E\Delta^T(\zeta)\right)\widetilde{M}^T \leq 0$$

$$\Rightarrow -Sym\{\widetilde{M}E^{-1}\widetilde{N}(rX + \bar{r}\bar{X})\} - \rho^{-2}\widetilde{M}E^{-1}\widetilde{M}^T - \rho^2(rX + \bar{r}\bar{X})^T \widetilde{N}^T E^{-1}\widetilde{N}(rX + \bar{r}\bar{X}) + Sym\{\widetilde{M}\Delta(\zeta)\widetilde{N}(rX + \bar{r}\bar{X})\} \leq 0 \quad (21)$$

According to (21), the inequality (19) holds if

$$Sym\{A_{0Cl}(rX + \bar{r}\bar{X}) + \widetilde{M}E^{-1}\widetilde{N}(rX + \bar{r}\bar{X})\} + \rho^{-2}\widetilde{M}E^{-1}\widetilde{M}^T + \rho^2(rX + \bar{r}\bar{X})^T \widetilde{N}^T E^{-1}\widetilde{N}(rX + \bar{r}\bar{X}) < 0 \quad (22)$$

Inequality (22) is equivalent to existence of positive constants $\rho$ and $\gamma$ such that

$$Sym\{A_{0Cl}(rX + \bar{r}\bar{X}) + \widetilde{M}E^{-1}\widetilde{N}(rX + \bar{r}\bar{X})\} + \rho^{-2}\widetilde{M}(E^{-1} + \gamma^{-1} I)\widetilde{M}^T$$

$$+\rho^2(rX+\bar{r}\bar{X})^T\tilde{N}^TE^{-1}\tilde{N}(rX+\bar{r}\bar{X})<0 \quad (23)$$

Inequality (23) can be rewritten as follows

$$Sym\{A_{0Cl}(rX+\bar{r}\bar{X})\}$$
$$+[\rho^{-1}\tilde{M}\ \ \rho(rX+\bar{r}\bar{X})^T\tilde{N}^T]\begin{bmatrix}E^{-1}+\gamma^{-1}I & E^{-1}\\ E^{-1} & E^{-1}\end{bmatrix}\begin{bmatrix}\rho^{-1}M\tilde{M}^T\\ \rho\tilde{N}(rX+\bar{r}\bar{X})\end{bmatrix}$$
$$=Sym\{A_{0Cl}(rX+\bar{r}\bar{X})\}$$
$$+[\rho^{-1}\tilde{M}\ \ \rho(rX+\bar{r}\bar{X})^T\tilde{N}^T]\begin{bmatrix}\gamma I & -\gamma I\\ -\gamma I & E+\gamma I\end{bmatrix}^{-1}\begin{bmatrix}\rho^{-1}M\tilde{M}^T\\ \rho\tilde{N}(rX+\bar{r}\bar{X})\end{bmatrix}$$
$$<0 \quad (24)$$

which can be rewritten as follows using Schur complement [28],

$$\begin{bmatrix}Sym\{A_{0Cl}(rX+\bar{r}\bar{X})\} & \rho^{-1}\tilde{M} & \rho(rX+\bar{r}\bar{X})^T\tilde{N}^T\\ \bullet & -\gamma I & \gamma I\\ \bullet & \bullet & -(E+\gamma I)\end{bmatrix}<0 \quad (25)$$

One can get the following inequality through pre- and post-multiplying (25) by $diag(\rho I, I, I)$

$$\begin{bmatrix}Sym\{A_{0Cl}(r(\rho^2 X)+\overline{r(\rho^2 X)})\} & \rho^{-1}\tilde{M} & \rho^2(r(\rho^2 X)+\overline{r(\rho^2 X)})^T\tilde{N}^T\\ \bullet & -\gamma I & \gamma I\\ \bullet & \bullet & -(E+\gamma I)\end{bmatrix}$$
$$<0 \quad (26)$$

Defining $P=\rho^2 X$ which is assumed to be in the form of (14), the inequality (26) can be written as follows

$$\begin{bmatrix}\Pi'_{11} & \tilde{M} & \Pi'_{13}\\ \bullet & -\gamma I & \gamma I\\ \bullet & \bullet & -(E+\gamma I)\end{bmatrix}<0 \quad (27)$$

where

$$\Pi'_{11}=\begin{bmatrix}p'_{11} & p'_{12}\\ p'_{21} & p'_{22}\end{bmatrix}, \Pi'_{13}=\begin{bmatrix}q'_1\\ q'_2\end{bmatrix}$$

$$p'_{11}=A(rP_S+\bar{r}\bar{P}_S)+(rP_S+\bar{r}\bar{P}_S)^T A^T+$$
$$BD_C C(rP_S+\bar{r}\bar{P}_S)+(rP_S+\bar{r}\bar{P}_S)^T C^T D_C^T B^T,$$
$$p'_{12}=BC_C(rP_C+\bar{r}\bar{P}_C)+(rP_C+\bar{r}\bar{P}_C)^T C_C^T B^T,$$
$$p'_{21}=B_C C(rP_S+\bar{r}\bar{P}_S)+(rP_S+\bar{r}\bar{P}_S)^T C_C^T B^T,$$
$$p'_{22}=A_C(rP_C+\bar{r}\bar{P}_C)+(rP_C+\bar{r}\bar{P}_C)^T A_C^T,$$
$$q'_1=(rP_S+\bar{r}\bar{P}_S)^T N_1^T+(rP_S+\bar{r}\bar{P}_S)^T C^T D_C^T N_2^T, q'_2$$
$$=(rP_C+\bar{r}\bar{P}_C)^T C_C^T N_2^T. \quad (28)$$

However, the matrix inequality (27) is not linear due to several multiplications of variables. Therefore, by linearizing change of variables as

$$T_1=A_C(rP_C+\bar{r}\bar{P}_C), T_2=B_C C(rP_S+\bar{r}\bar{P}_S),$$
$$T_3=C_C(rP_C+\bar{r}\bar{P}_C), T_4=D_C C(rP_S+\bar{r}\bar{P}_S), \quad (29)$$

equations in (28) can be rewritten as

$$p'_{11}=A(rP_S+\bar{r}\bar{P}_S)+(rP_S+\bar{r}\bar{P}_S)^T A^T+BT_4+T_4^T B^T,$$
$$p'_{12}=BT_3+T_2^T,\ \ p'_{21}=T_2+T_3^T B^T,$$
$$p'_{22}=T_1+T_1^T, q'_1=(rP_S+\bar{r}\bar{P}_S)^T N_1^T+T_4^T N_2^T,$$
$$q'_2=T_3^T N_2^T. \quad (30)$$

by which the LMI in (27) is equivalent to the linear matrix inequality in (15). ∎

**Theorem 2**: Considering closed-loop system in (12) with $1\leq\alpha<2$, and a positive definite symmetric matrix $P=P^T$ in the form of (14) with $P_S\in\mathcal{R}^{n\times n}$ and $P_C\in\mathcal{R}^{n_C\times n_C}$ and a real scalar constant $\gamma>0$ alongside with matrices $T_i, i=1,\ldots,4$ exist such that the following LMI constrain become feasible

$$\begin{bmatrix}\hat{\Pi}_{11} & \Theta\otimes\tilde{M} & \hat{\Pi}_{13}\\ \bullet & -\gamma I & \gamma I\\ \bullet & \bullet & -Sym(I_2\otimes J)-\gamma I\end{bmatrix}<0, \quad (31)$$

in which

$$\hat{\Pi}_{11}=\begin{bmatrix}\hat{p}_{11} & \hat{p}_{12}\\ \hat{p}_{21} & \hat{p}_{22}\end{bmatrix}, \hat{\Pi}_{13}=\begin{bmatrix}\hat{q}_1 & 0\\ 0 & \hat{q}_2\end{bmatrix},$$

$$\hat{p}_{11}=\hat{p}_{22}$$
$$=\begin{bmatrix}AP_S+P_S A^T+BT_4+T_4^T B^T & BT_3+T_2^T\\ T_2+T_3^T B^T & T_1+T_1^T\end{bmatrix}sin\theta,$$

$$\hat{p}_{12}=-\hat{p}_{21}$$
$$=\begin{bmatrix}P_S A^T-AP_S+T_4^T B^T-BT_4 & T_2^T-BT_3\\ T_3^T B^T-T_2 & T_1^T-T_1\end{bmatrix}cos\theta,$$

$$\hat{q}_1=\hat{q}_2=\begin{bmatrix}P_S N_1^T+T_4^T N_2^T\\ T_3^T N_2^T\end{bmatrix} \quad (32)$$

where $\theta=\pi-\alpha\pi/2$ and $\Theta$ is defined in (8) then, the dynamic output feedback controller parameters of

$$A_C=T_1 P_C^{-1}, B_C=T_2 P_S^{-1} C^\dagger,$$
$$C_C=T_3 P_C^{-1}, D_C=T_4 P_S^{-1} C^\dagger, \quad (33)$$

make the closed-loop system in (12) asymptotically stable.

**Proof**: It follows from Lemma 2 that the uncertain fractional-order closed-loop system (12) with $1<\alpha<2$ and $A_{Cl}=A_{0Cl}+A_{\Delta Cl}$, defined in (18), is asymptotically stable if there exists a positive definite matrix $X=X^T$, $X\in\mathcal{R}^{(n+n_C)\times(n+n_C)}$ such that

$$Sym\{\Theta\otimes(A_{Cl}X)\}<0\Leftrightarrow Sym\{\Theta\otimes\bigl((A_{0Cl}+A_{\Delta Cl})X\bigr)\}$$
$$<0$$

$$\Leftrightarrow Sym\{\Theta\otimes(A_{0Cl}X)\}+Sym\{\Theta\otimes\bigl(\tilde{M}\Delta(\zeta)\tilde{N}X\bigr)\}<0$$

$$\Leftrightarrow\Sigma+Sym\{\hat{M}\hat{\Delta}(\zeta)\hat{N}\hat{X}\}<0, \quad (34)$$

where $\Sigma=Sym\{\Theta\otimes(A_{0Cl}X)\}, \hat{M}=\Theta\otimes\tilde{M}, \hat{N}=I_2\otimes\tilde{N}, \hat{F}(\zeta)=I_2\otimes F(\zeta), \hat{J}=I_2\otimes J, \hat{\Delta}(\zeta)=I_2\otimes\Delta(\zeta)=\hat{F}(\zeta)[I-\hat{J}\hat{F}(\zeta)]^{-1}, \hat{X}=I_2\otimes X$. Introducing

$$\hat{E}=Sym(\hat{J}), \hat{F}=\hat{E}^{-\frac{1}{2}}(\rho^{-1}\hat{M}^T+\rho\hat{N}X)$$
$$-\hat{E}^{\frac{1}{2}}\Delta^T(\zeta)\rho^{-1}\hat{M}^T, \quad (35)$$

it follows from Lemma 3 that $Sym\{\hat{\Delta}(\zeta)\}-\hat{\Delta}(\zeta)\hat{E}\hat{\Delta}^T(\zeta)>0$. Therefore, the following inequality holds for any $\rho>0$

$$-\hat{F}^T\hat{F} \leq 0 \Leftrightarrow -\left(\hat{E}^{-\frac{1}{2}}(\rho^{-1}\widehat{M}^T + \rho\widehat{N}X)\right.$$
$$\left. - \hat{E}^{\frac{1}{2}}\Delta^T(\zeta)\rho^{-1}\widehat{M}^T\right)^T\left(\hat{E}^{-\frac{1}{2}}(\rho^{-1}\widehat{M}^T\right.$$
$$\left. + \rho\widehat{N}X) - \hat{E}^{\frac{1}{2}}\Delta^T(\zeta)\rho^{-1}\widehat{M}^T\right) < 0 \Leftrightarrow$$
$$-Sym\{\widehat{M}\hat{E}^{-1}\widehat{N}X\} - \rho^{-2}\widehat{M}\hat{E}^{-1}\widehat{M}^T - \rho^2 X^T\widehat{N}\hat{E}^{-1}\widehat{N}X$$
$$+ Sym\{\widehat{M}\hat{\Delta}(\zeta)\widehat{N}X\}$$
$$+ \rho^{-2}\widehat{M}\left(Sym\{\hat{\Delta}(\zeta)\} - \hat{\Delta}(\zeta)\hat{E}\hat{\Delta}^T(\zeta)\right)\widehat{M}^T$$
$$< 0$$
$$\Rightarrow -Sym\{\widehat{M}\hat{E}^{-1}\widehat{N}X\} - \rho^{-2}\widehat{M}\hat{E}^{-1}\widehat{M}^T - \rho^2 X^T\widehat{N}^T\hat{E}^{-1}\widehat{N}X$$
$$+ Sym\{\widehat{M}\,\hat{\Delta}(\zeta)\widehat{N}X\}$$
$$\leq 0 \quad (36)$$

According to (36), the inequality (35) holds if

$$\Sigma + Sym\{\widehat{M}\hat{E}^{-1}\widehat{N}X\} + \rho^{-2}\widehat{M}\hat{E}^{-1}\widehat{M}^T + \rho^2 X^T\widehat{N}^T\hat{E}^{-1}\widehat{N}X < 0 \quad (37)$$

which is equivalent to existence of positive constants $\rho$ and $\gamma$ such that

$$\Sigma + Sym\{\widehat{M}\hat{E}^{-1}\widehat{N}X\} + \rho^{-2}\widehat{M}(\hat{E}^{-1} + \gamma I)\widehat{M}^T + \rho^2 X^T\widehat{N}^T\hat{E}^{-1}\widehat{N}X < 0 \quad (38)$$

Inequality (38) can be rewritten as follows

$$\Sigma + [\rho^{-1}\widehat{M} \quad \rho\hat{E}^{-1}X^T\widehat{N}^T]\begin{bmatrix}\hat{E}^{-1}+\gamma^{-1}I & \hat{E}^{-1} \\ \hat{E}^{-1} & \hat{E}^{-1}\end{bmatrix}\begin{bmatrix}\rho^{-1}\widehat{M}^T \\ \rho\widehat{N}X\end{bmatrix}$$
$$= \Sigma + [\rho^{-1}\widehat{M} \quad \rho\hat{E}^{-1}X^T\widehat{N}^T]\begin{bmatrix}\gamma I & -\gamma I \\ -\gamma I & \hat{E}+\gamma I\end{bmatrix}\begin{bmatrix}\rho^{-1}\widehat{M}^T \\ \rho\widehat{N}X\end{bmatrix}$$
$$< 0 \quad (39)$$

which can be, by Schur complement [28], rewritten as

$$\begin{bmatrix}\Sigma & \rho^{-1}\widehat{M} & \rho X^T\widehat{N}^T \\ \bullet & -\gamma I & \gamma I \\ \bullet & \bullet & -(\hat{E}+\gamma I)\end{bmatrix} < 0 \quad (40)$$

One can get the following inequality through pre- and post-multiplying (40) by $diag(\rho I, I, I)$

$$\begin{bmatrix}\rho^2\Sigma & \widehat{M} & \rho^2 X^T\widehat{N}^T \\ \bullet & -\gamma I & \gamma I \\ \bullet & \bullet & -(\hat{E}+\gamma I)\end{bmatrix} < 0 \quad (41)$$

Defining $P = \rho^2\hat{X}$ in which $X$ is assumed to be in the form of (14), the inequality (41) can be written as follows

$$\begin{bmatrix}\widehat{\Pi}'_{11} & \widehat{M} & \widehat{\Pi}'_{13} \\ \bullet & -\gamma I & \gamma I \\ \bullet & \bullet & -(\hat{E}+\gamma I)\end{bmatrix} < 0 \quad (42)$$

where

$$\widehat{\Pi}'_{11} = \begin{bmatrix}\hat{p}'_{11} & \hat{p}'_{12} \\ \hat{p}'_{21} & \hat{p}'_{22}\end{bmatrix}, \widehat{\Pi}'_{13} = \begin{bmatrix}\hat{q}' & 0 \\ 0 & \hat{q}'\end{bmatrix}$$

$$\hat{p}'_{11} = \hat{p}'_{22}$$
$$= \begin{bmatrix}AP_S + P_S A^T + BD_C CP_S + P_S C^T D_C^T B^T & BC_C P_C + P_S C^T B_C^T \\ B_C CP_S + P_C C_C^T B^T & A_C P_C + P_C A_C^T\end{bmatrix}\sin\theta$$

$$\hat{p}'_{12} = -\hat{p}'_{21}$$
$$= \begin{bmatrix}P_S A^T - AP_S + P_S C^T D_C^T B^T - BD_C CP_S & P_S C^T B_C^T - BC_C P_C \\ P_C C_C^T B^T - B_C CP_S & P_C A_C^T - A_C P_C\end{bmatrix}\cos\theta$$

$$\hat{q}' = \begin{bmatrix}P_S N_1^T + P_S C^T D_C^T N_2^T \\ P_C C_C^T N_2^T\end{bmatrix} \quad (43)$$

The matrix inequality (42) is not linear because of various multiplications of variables. Accordingly, by changing variables as follows

$$T_1 = A_C P_C, \quad T_2 = B_C CP_S,$$
$$T_3 = C_C P_C, \quad T_4 = D_C CP_S, \quad (44)$$

equations in (43) can be rewritten as

$$\hat{p}'_{11} = \hat{p}'_{22} = \begin{bmatrix}AP_S + P_S A^T + BT_4 + T_4^T B^T & BT_3 + T_2^T \\ T_2 + T_3^T B^T & T_1 + T_1^T\end{bmatrix}\sin\theta$$

$$\hat{p}'_{12} = -\hat{p}'_{21} = \begin{bmatrix}P_S A^T - AP_S + T_4^T B^T - BT_4 & T_2^T - BT_3 \\ T_3^T B^T - T_2 & T_1 - T_1^T\end{bmatrix}\cos\theta$$

$$\hat{q}' = \begin{bmatrix}P_S N_1^T + T_4^T N_2^T \\ T_3^T N_2^T\end{bmatrix} \quad (45)$$

by which the LMI in (42) is equivalent to the linear matrix inequality in (31). ∎

**Remark 2**: By solving the proposed LMIs in Theorem 1 and Theorem 2 for special case of $n_C = 0$, the static output feedback controllers are obtained for $0 < \alpha < 1$ and $1 \leq \alpha < 2$ respectively.

**Remark 3**: in Theorem 1 and 2, control parameters are obtained by finding a proper block-diagonal positive definite matrix $P$, which can result in conservatism of the problem.

**Corollary 1**: Although Theorem 1 and Theorem 2 are allocated to robust stabilization of uncertain FO-LTI systems of form (1), the proposed method can be easily used for the case of certain systems by solving the LMI constraints $\Pi_{11} < 0$ and $\widehat{\Pi}_{11} < 0$ in these theorems, respectively.

**Proof:** The proof is straightforward by assuming $M = \mathbf{0}$ in proof procedure of Theorem 1 and Theorem 2.

## 4. Numerical examples

In this section, some numerical examples are given to show the applicability of the derived results. Several helpful LMI Parsers and Solvers can be used to assess the feasibility of the proposed constraints in order to obtain the design parameters. In this paper, we use YALMIP parser [29] and SeDuMi [30] solver in Matlab tool [31].

### 4.1. Example 1

Consider the dynamic output feedback stabilization problem of the uncertain fractional-order system (1) with the following parameters, which is also available in [10].

$$A = \begin{bmatrix}2.5 & 10 & -5 \\ -15 & 7.5 & -5 \\ 10 & 10 & -5\end{bmatrix}, M = \begin{bmatrix}0.5 & 1 & 0 \\ -0.4 & 0.2 & 0 \\ 0.1 & -0.1 & -0.6\end{bmatrix},$$

$$N_1 = \begin{bmatrix} 0.5 & 1.5 & 2 \\ 0 & 0.5 & 2.5 \\ 0 & 0 & 2.5 \end{bmatrix}, \quad B = \begin{bmatrix} 1 \\ 2 \\ 1 \end{bmatrix}, N_2 = \begin{bmatrix} 1 \\ -0.5 \\ 0.5 \end{bmatrix},$$
$$C = \begin{bmatrix} 1 & 1 & 1 \\ 0 & -2 & -2 \end{bmatrix}, J = I_3, \alpha = 0.8. \tag{46}$$

According to Theorem 1, it can be concluded that the uncertain fractional-order system (1) with the parameters in (46) is asymptotically stabilizable using the obtained dynamic output feedback controllers of arbitrary orders, tabulated in Table 1. Fig. 1 illustrates the time response of the resulted uncertain closed-loop FO-LTI system of form (12) via obtained controllers with $n_c = 1$ and 3, and the static controller in [10]. It is shown in Fig. 1 that all the states asymptotically converge to zero. It can be concluded from Fig. 1 that the oscillation and settling time of the response of the uncertain system via obtained dynamic output feedback controller with $n_c = 1$ and 3, are better than the that of the static controller proposed in [10] for the same problem.

**Table 1** Controller parameters obtained by Theorem 1 for Example 1.

| $n_c$ | $A_c$ | $B_c$ | $C_c$ | $D_c$ |
|---|---|---|---|---|
| 0 | 0 | 0 | 0 | $\begin{bmatrix} 0.2 \\ 0.7 \end{bmatrix}^T$ |
| 1 | $-45.4$ | $[-1.1 \; -0.8]$ | 1.1 | $\begin{bmatrix} -6.4 \\ -1.8 \end{bmatrix}^T$ |
| 2 | $\begin{bmatrix} -28.7 & 0.0 \\ 0.5 & -28.4 \end{bmatrix}$ | $\begin{bmatrix} -0.1 & 0.0 \\ -1.3 & -0.6 \end{bmatrix}$ | $\begin{bmatrix} 0.8 \\ 0.9 \end{bmatrix}^T$ | $\begin{bmatrix} -5.8 \\ -1.9 \end{bmatrix}^T$ |
| 3 | $\begin{bmatrix} -24.9 & 0.0 & 0.0 \\ 0.6 & -25.5 & -0.0 \\ 0.6 & -0.0 & -25.5 \end{bmatrix}$ | $\begin{bmatrix} 2.1 & 1.4 \\ 1.9 & 0.9 \\ 1.6 & 0.7 \end{bmatrix}$ | $\begin{bmatrix} 1.0 \\ 0.1 \\ 0.1 \end{bmatrix}^T$ | $\begin{bmatrix} -3.6 \\ -0.7 \end{bmatrix}^T$ |

Moreover, the obtained controllers' efforts are depicted in Fig. 2 in order to have a more precise comparison. It can be seen that, with almost the same amount of control effort, the obtained dynamic output feedback controllers, even with a low order of $n_c = 1$, have more appropriate stabilizing actions compared to static one proposed in [10].

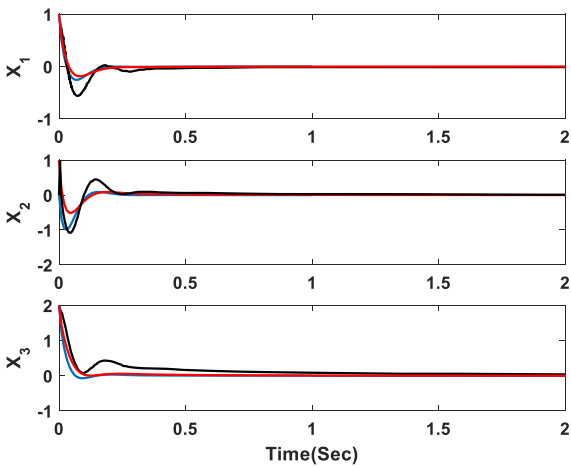

**Fig. 1.** *The time response of the closed-loop system in Example 1 via obtained output feedback controllers with $\boldsymbol{n_c = 1}$ (blue) $\boldsymbol{n_c = 3}$ (red) and the static controller obtained in [10] (black).*

In the following, robustness of the proposed method has been validated through numerical simulation results of 50 random systems with positive real uncertainties defined in (46) which have been stabilized via obtained dynamic output controller with $n_c = 2$ whose output time response are depicted in Fig. 3.

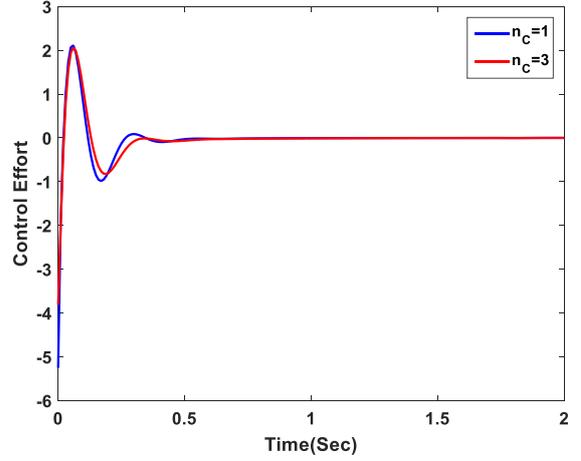

**Fig. 2.** *The control effort of dynamic output feedback controllers obtained by Theorem 1.*

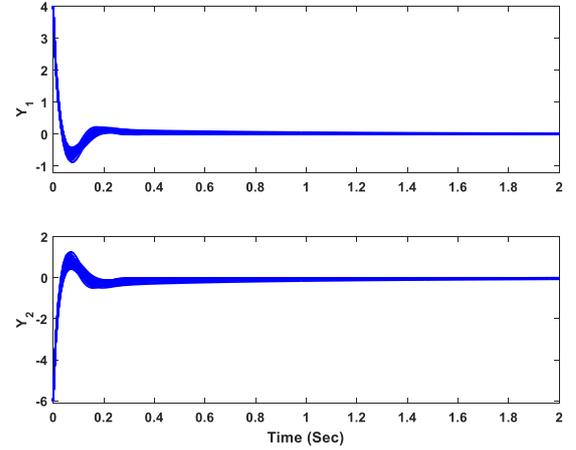

**Fig. 3.** *Output time responses of 50 random systems with positive real uncertainties defined in Example 1 obtained via dynamic output controller with $\boldsymbol{n_c = 2}$.*

### 4.2. Example 2

Dynamic output feedback stabilization problem of the uncertain fractional-order system (1) with the following parameters is considered, which is also available in [10].

$$A = \begin{bmatrix} 0 & 4 & 6 & 6 \\ 0 & -4 & 2 & 0 \\ -2 & -2 & -6 & 2 \\ 10 & 10 & 4 & 0 \end{bmatrix}, M = \begin{bmatrix} 0.2 & 0.2 & 0 & 0 \\ 0 & 0.2 & 0.2 & 0 \\ 0 & 0 & 0.2 & 0.2 \\ 0 & 0 & 0 & 0.2 \end{bmatrix},$$
$$N_1 = \begin{bmatrix} 0.1 & 0.3 & 0.4 & 0 \\ 0 & 0.1 & 0.5 & 0 \\ 0 & 0 & 0.5 & 0 \\ 0 & 0 & 0 & 1 \end{bmatrix}, B = \begin{bmatrix} 1 \\ 2 \\ 1 \\ 0 \end{bmatrix}, N_2 = \begin{bmatrix} 0.2 \\ -0.1 \\ 0.1 \\ 0.1 \end{bmatrix},$$
$$C = \begin{bmatrix} 5 & 2 & 3 & 10 \\ -20 & 15 & 0 & 20 \end{bmatrix}, J = I_4, \alpha = 1.2. \tag{47}$$

According to Theorem 2, we can conclude that the uncertain fractional-order system (1), with the parameters in (47), is asymptotically stabilizable using the obtained dynamic output feedback controllers of arbitrary orders, in the form of (11), proposed in Table 2. The time response of the resulted uncertain closed-loop FO-LTI system of form (12) via obtained controllers with $n_c = 1$ and 4, and the proposed static controller in [10], are illustrated in Fig. 4 in which, all the states asymptotically converge to zero. It can be concluded from Fig. 3 that the oscillation and settling time of the response of the uncertain system via obtained dynamic output feedback controller with $n_c = 1$ and 4, are far better than the that of the static controller proposed in [10] for the same problem.

**Table 2** Controller parameters obtained by Theorem 2 for Example 2.

| $n_c$ | $A_c$ | $B_c$ | $C_c$ | $D_c$ |
|---|---|---|---|---|
| 0 | 0 | 0 | 0 | $\begin{bmatrix}-2\\0.0\end{bmatrix}^T$ |
| 1 | $-309.6$ | $[25.4\ \ -1.3]$ | $16.5$ | $\begin{bmatrix}-3.3\\0.2\end{bmatrix}^T$ |
| 2 | $\begin{bmatrix}-134.7 & -55.2\\-55.3 & -134.5\end{bmatrix}$ | $\begin{bmatrix}34.8 & -1.6\\34.4 & -1.5\end{bmatrix}$ | $\begin{bmatrix}10.2\\10.6\end{bmatrix}^T$ | $\begin{bmatrix}-5.3\\0.3\end{bmatrix}^T$ |
| 3 | $\begin{bmatrix}-91.8 & -39.3 & -39.8\\-39.3 & 91.7 & -39.6\\-39.7 & -39.3 & -91.7\end{bmatrix}$ | $\begin{bmatrix}55.4 & -2.8\\55.1 & -2.3\\55.9 & -2.5\end{bmatrix}$ | $\begin{bmatrix}6.0\\5.9\\5.8\end{bmatrix}^T$ | $\begin{bmatrix}-7.6\\0.4\end{bmatrix}^T$ |
| 4 | $\begin{bmatrix}-62.5 & -27.1 & -27.1 & -27.2\\-27.1 & -62.2 & -27.6 & -27.1\\-27.7 & -27.2 & -62.1 & -27.2\\-27.2 & -27.7 & -27.7 & -62.5\end{bmatrix}$ | $\begin{bmatrix}65.9 & -3.7\\66.0 & -3.3\\65.1 & -3.4\\65.9 & -3.4\end{bmatrix}$ | $\begin{bmatrix}3.4\\3.4\\3.3\\3.5\end{bmatrix}^T$ | $\begin{bmatrix}-7.5\\0.4\end{bmatrix}^T$ |

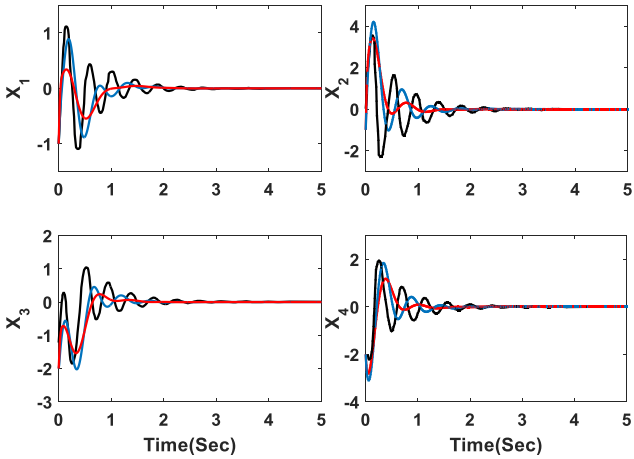

**Fig. 4.** *The time response of the closed-loop system in Example 2 via obtained output feedback controllers with $n_c = 1$ (blue) $n_c = 4$ (red) and the static controller obtained in [10] (black).*

Moreover, the obtained controllers' efforts are depicted in Fig. 5, where with almost the same amount of control effort, the obtained dynamic output feedback controllers, even with a low order of $n_c = 1$, have more appropriate stabilizing actions compared to static one proposed in [10].

Besides, robustness of the proposed method has been validated through numerical simulation results of 50 random systems with positive real uncertainties defined in (47) which have been stabilized via obtained dynamic output controller with $n_c = 2$ whose output time response are depicted in Fig. 6.

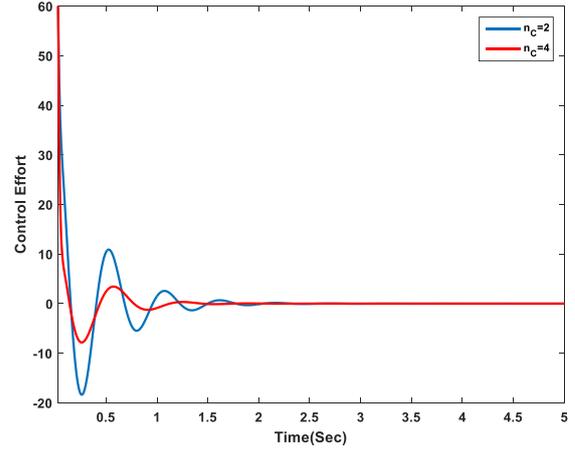

**Fig. 5.** *The control effort of dynamic output feedback controllers obtained by Theorem 2.*

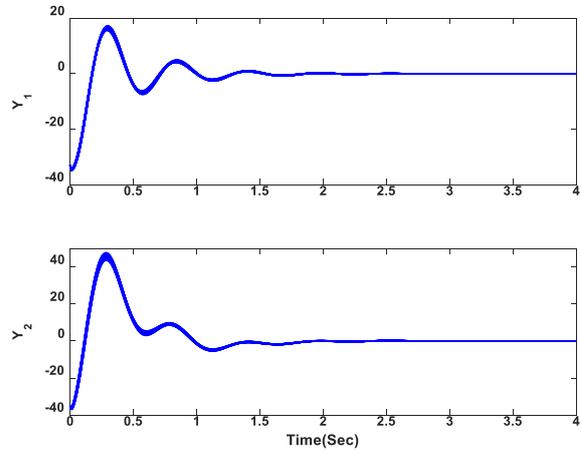

**Fig. 6.** *Output time responses of 50 random systems with positive real uncertainties defined in Example 2 obtained via dynamic output controller with $\mathbf{n_c = 2}$.*

### 4.3. Example 3

Consider the dynamic output feedback stabilization problem of the uncertain fractional-order multi-input and multi-output system (1) with the following parameters

$$A = \begin{bmatrix}2.5 & 10 & -5\\-15 & 7.5 & -5\\10 & 10 & -5\end{bmatrix}, M = \begin{bmatrix}0.5 & 1 & 0\\-0.4 & 0.2 & 0\\0.1 & -0.1 & -0.6\end{bmatrix},$$
$$N_1 = \begin{bmatrix}0.5 & 1.5 & 2\\0 & 0.5 & 2.5\\0 & 0 & 2.5\end{bmatrix}, B = \begin{bmatrix}1.3 & 1\\2 & 0\\0.9 & -0.8\end{bmatrix},$$
$$N_2 = \begin{bmatrix}1 & 0.5\\-0.5 & -1\\0.5 & 0.5\end{bmatrix}, C = \begin{bmatrix}1 & 1 & 0\\0 & 0 & -1\\1 & 1 & -1\end{bmatrix}, J = I_3, \alpha = 0.9. \quad (48)$$

It should be noted that no static output and Observer-based output feedback stabilization controller can be designed using the method proposed in [10], because of rank constraint on output matrix $C$. Nevertheless, using our proposed method static output feedback controller and dynamic output feedback controllers with predetermined orders can be easily obtained. According to Theorem 1, it can be concluded that the uncertain fractional-order system (1) with the parameters in (48) is asymptotically stabilizable using the obtained dynamic output feedback controllers of arbitrary orders, tabulated in Table 3. Fig. 7 illustrates the time response of the resulted uncertain closed-loop FO-LTI system of form (12) via obtained controllers with $n_c = 0,1$ and 3, where all the states asymptotically converge to zero. However, the state trajectories of the system via obtained static controller for $n_c = 0$ is very oscillatory.

Table **3** Controller parameters obtained by Theorem 1 for Example 3.

| $n_c$ | $A_c$ | $B_c$ | $C_c$ | $D_c$ |
|---|---|---|---|---|
| 0 | 0 | 0 | 0 | $\begin{bmatrix} -1.49 & 1.11 \\ 0.14 & -0.43 \\ -1.35 & 0.68 \end{bmatrix}^T$ |
| 1 | $-45.44$ | $\begin{bmatrix} -3.92 \\ 3.70 \\ -0.22 \end{bmatrix}^T$ | $\begin{bmatrix} 4.105 \\ -2.80 \end{bmatrix}$ | $\begin{bmatrix} -8.76 & 6.17 \\ 6.78 & -5.86 \\ -1.97 & 0.31 \end{bmatrix}^T$ |
| 2 | $\begin{bmatrix} -28.15 & -0.00 \\ 0.84 & -28.35 \end{bmatrix}$ | $\begin{bmatrix} -1.85 & -2.16 \\ 1.58 & 1.96 \\ -0.27 & -0.20 \end{bmatrix}^T$ | $\begin{bmatrix} 1.74 & 1.87 \\ -1.54 & -1.68 \end{bmatrix}$ | $\begin{bmatrix} -1.49 & 1.11 \\ 0.14 & -0.43 \\ -1.35 & 0.68 \end{bmatrix}^T$ |
| 3 | $\begin{bmatrix} -24.37 & 0.00 & -0.00 \\ 1.10 & -23.73 & 0.00 \\ 1.47 & 2.21 & -23.48 \end{bmatrix}$ | $\begin{bmatrix} -1.58 & 1.65 & 0.07 \\ -2.82 & 2.64 & -0.17 \\ -3.63 & 3.32 & -0.31 \end{bmatrix}^T$ | $\begin{bmatrix} 1.83 & -1.45 \\ 3.02 & -2.11 \\ 3.53 & -2.18 \end{bmatrix}^T$ | $\begin{bmatrix} -5.61 & 2.82 \\ 3.88 & -2.88 \\ -1.73 & -0.06 \end{bmatrix}^T$ |

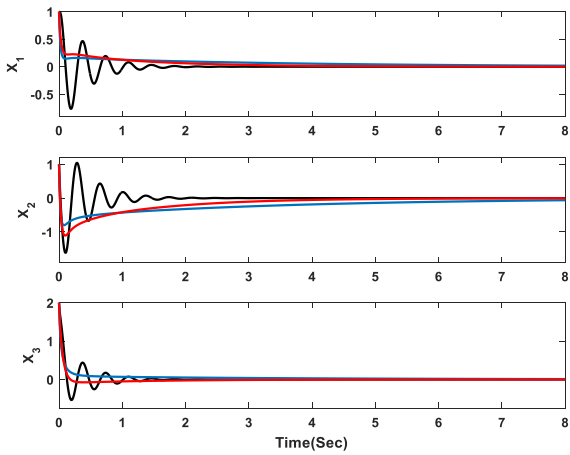

**Fig. 7.** *The time response of the closed-loop system in Example 3 via obtained output feedback controllers with $\mathbf{n_c = 0}$ (black), $\mathbf{n_c = 1}$ (blue) $\mathbf{n_c = 3}$ (red).*

In addition, the obtained controllers' efforts are depicted in Fig. 8 in order to have a more precise comparison, in which the control effort of the static controller is also very oscillatory. On the other hand, with almost the same amount of control effort, the obtained dynamic output feedback controllers, even with a low order of $n_c = 1$, have more appropriate stabilizing actions compared to static one.

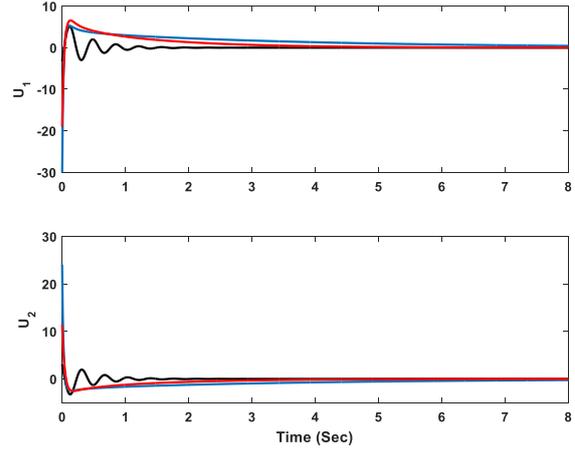

**Fig. 8.** *The control effort of static and dynamic output feedback controllers with $\mathbf{n_c = 0}$ (black), $\mathbf{n_c = 1}$ (blue) $\mathbf{n_c = 3}$ (red), obtained by Theorem 1.*

In the following, robustness of the proposed method has been demonstrated through numerical simulation results of 50 random systems with positive real uncertainties defined in (48) which have been stabilized via obtained dynamic output controller with $n_c = 2$ whose output time response are depicted in Fig. 9.

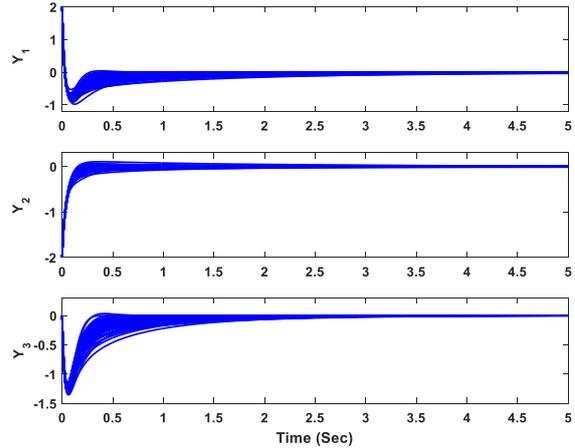

**Fig. 9.** *Output time responses of 50 random systems with positive real uncertainties defined in Example 3 obtained via dynamic output controller with $\mathbf{n_c = 2}$.*

## 5. Conclusion

In this paper a fixed-order dynamic output feedback controller design method based on LMI approach is proposed to robustly stabilize fractional-order linear systems with positive real uncertainty. Sufficient conditions are obtained for designing a stabilizing controller with a predetermined order, which can be chosen to be as low as possible for simpler implementation. Indeed by using proposed method, one can benefit from dynamic output feedback controller advantages with orders lower than the system order. The LMI-based approach of developing robust stabilizing control is preserved in spite of the complexity of assuming the most complete model of linear controller, including direct feedthrough parameter. Moreover, no limiting constraint on the state space matrices of the uncertain systems is considered in our controller design procedure. Eventually, some numerical examples are presented to show the effectiveness of the proposed controller design method. An area of future work is to find a method to obtain the controller parameters using a full positive definite matrix to avoid the potential conservatism of the proposed procedure. Moreover, comprehensiveness of the proposed method can be enhanced considering different types of uncertainties, such as polytopic and interval uncertainty.